\documentclass[a4paper,10pt]{amsart}

\usepackage[latin1]{inputenc}
\usepackage{amsmath}
\usepackage{amsfonts}
\usepackage{amssymb}
\usepackage[T1]{fontenc}
\usepackage[english]{babel}
\usepackage{dsfont}
\usepackage[all]{xy}
\usepackage{stmaryrd}

\newcommand{\zb}{\mathbb{Z}}

\newcommand{\lb}{\mathbb{L}}

\newcommand{\al}{\alpha}

\newcommand{\fl}{\rightarrow}
\newcommand{\gfl}{\longrightarrow}
\newcommand{\infl}{\rightarrowtail}
\newcommand{\defl}{\twoheadrightarrow}

\newcommand{\dr}{\ar@{->}[r]}
\newcommand{\dri}{\ar@{>->}[r]}
\newcommand{\drp}{\ar@{-->}[r]}
\newcommand{\dre}{\ar@{->>}[r]}
\newcommand{\dreg}{\ar@{=}[r]}
\newcommand{\drm}{\ar@{^{(}->}[r]}
\newcommand{\ddr}{\ar@{->}[rr]}
\newcommand{\ddre}{\ar@{->>}[rr]}
\newcommand{\ddreg}{\ar@{=}[rr]}
\newcommand{\ha}{\ar@{->}[u]}
\newcommand{\hae}{\ar@{->>}[u]}
\newcommand{\hap}{\ar@{-->}[u]}
\newcommand{\ham}{\ar@{^{(}->}[u]}
\newcommand{\hham}{\ar@{^{(}->}[uu]}
\newcommand{\hag}{\ar@{->}[ul]}
\newcommand{\hagm}{\ar@{^{(}->}[ul]}
\newcommand{\hagp}{\ar@{-->}[ul]}
\newcommand{\hdr}{\ar@{->}[ur]}
\newcommand{\hdrm}{\ar@{^{(}->}[ur]}
\newcommand{\hdri}{\ar@{>->}[ur]}
\newcommand{\hdre}{\ar@{->>}[ur]}
\newcommand{\bas}{\ar@{->}[d]}
\newcommand{\bbas}{\ar@{->}[dd]}
\newcommand{\basm}{\ar@{^{(}->}[d]}
\newcommand{\basp}{\ar@{-->}[d]}
\newcommand{\baseg}{\ar@{=}[d]}
\newcommand{\bbaseg}{\ar@{=}[dd]}
\newcommand{\bdr}{\ar@{->}[dr]}
\newcommand{\bdre}{\ar@{->>}[dr]}
\newcommand{\bbdr}{\ar@{->}[ddr]}
\newcommand{\bddr}{\ar@{->}[drr]}
\newcommand{\bg}{\ar@{->}[dl]}
\newcommand{\bgm}{\ar@{^{(}->}[dl]}
\newcommand{\bgp}{\ar@{-->}[dl]}
\newcommand{\bggp}{\ar@{-->}[dll]}
\newcommand{\bbgp}{\ar@{-->}[ddl]}

\newcommand{\ac}{\mathcal{A}}
\newcommand{\bc}{\mathcal{B}}
\newcommand{\cat}{\mathcal{C}}

\newcommand{\ec}{\mathcal{E}}
\newcommand{\hc}{\mathcal{H}}
\newcommand{\mc}{\mathcal{M}}
\newcommand{\pc}{\mathcal{P}}
\newcommand{\rc}{\mathcal{R}}
\newcommand{\tc}{\mathcal{T}}
\newcommand{\xc}{\mathcal{X}}
\newcommand{\es}{\underline{\mathcal{E}}}
\newcommand{\ms}{\underline{\mathcal{M}}}

\newcommand{\dba}{\mathcal{D}^\text{b}(\operatorname{mod} A)}
\newcommand{\eac}{\mathcal{H}^b_{\mathcal{E}-ac}\left(\mathcal{M}\right)}
\newcommand{\eace}{\mathcal{H}^b_{\mathcal{E}-ac}\left(\mathcal{E}\right)}
\newcommand{\hbm}{\mathcal{H}^b\left(\mathcal{M}\right)}
\newcommand{\hbms}{\mathcal{H}^b\left(\underline{\mathcal{M}}\right)}
\newcommand{\hbmp}{\mathcal{H}^b\left(\mathcal{M'}\right)}
\newcommand{\hbe}{\mathcal{H}^b\left(\mathcal{E}\right)}
\newcommand{\hbp}{\mathcal{H}^b\left(\mathcal{P}\right)}
\newcommand{\dbe}{\mathcal{D}^b\left(\mathcal{E}\right)}
\newcommand{\der}{\mathcal{D}\operatorname{Mod}\,\mathcal{M}}
\newcommand{\derp}{\mathcal{D}\operatorname{Mod}\,\mathcal{M'}}
\newcommand{\modm}{\operatorname{Mod}\mathcal{M}}
\newcommand{\modms}{\operatorname{mod}\underline{\mathcal{M}}}
\newcommand{\modmsp}{\operatorname{mod}\underline{\mathcal{M}}'}

\newcommand{\projm}{\operatorname{proj}\mc}

\newcommand{\projp}{\operatorname{proj}\pc}
\newcommand{\perm}{\operatorname{per}\mathcal{M}}
\newcommand{\permp}{\operatorname{per}\mathcal{M'}}
\newcommand{\perP}{\operatorname{per}\mathcal{P}}
\newcommand{\perms}{\operatorname{per}_{\underline{\mathcal{M}}}\mathcal{M}}
\newcommand{\permsp}{\operatorname{per}_{\underline{\mathcal{M}}'}\mathcal{M'}}

\newcommand{\susp}{\Sigma}
\newcommand{\susm}{\Sigma^{-1}}
\newcommand{\ke}{\operatorname{Ker}}
\newcommand{\im}{\operatorname{Im}}

\newcommand{\ext}{\operatorname{Ext}}
\newcommand{\homph}{\operatorname{Hom}}
\newcommand{\Hom}{\operatorname{\mathcal{H}om}}
\newcommand{\lf}{\mathbb{L}F}
\newcommand{\kzero}{\operatorname{K}_{0}}
\newcommand{\kred}{\overline{\operatorname{K}}_{0}(\cat_A)}
\newcommand{\coh}{\operatorname{H}}
\newcommand{\komodms}{\kzero(\modms)}
\newcommand{\komodmsp}{\kzero(\modmsp)}
\newcommand{\koprojmsurp}{\kzero(\projm)/\kzero(\projp)}
\newcommand{\koprojmpsurp}{\kzero(\projm')/\kzero(\projp)}

\newtheorem{theo}{Theorem}
\newtheorem{prop}[theo]{Proposition}
\newtheorem{lem}[theo]{Lemma}
\newtheorem{cor}[theo]{Corollary}

\begin{document}

\title[Grothendieck group and mutation rule]
{Grothendieck group and generalized mutation 
rule for 2-Calabi--Yau triangulated categories}
\author{Yann Palu}
\address{Universit{\'e} Paris 7 - Denis Diderot,
UMR 7586 du CNRS, case 7012, 2 place Jussieu,
75251 Paris Cedex 05, France.}

\email{
\begin{minipage}[t]{5cm}
palu@math.jussieu.fr
\end{minipage}
}

\begin{abstract}
 We compute the Grothendieck group
of certain 2-Calabi--Yau
triangulated categories appearing
naturally in the study of the link
between quiver representations
and Fomin--Zelevinsky's cluster
algebras. In this setup,
we also prove a generalization of
Fomin--Zelevinsky's mutation rule.
\end{abstract}

\maketitle

\section*{Introduction}

In their study~\cite{CK1}
of the connections between
cluster algebras
(see~\cite{Zel})
and
quiver representations,
P. Caldero and B. Keller
conjectured
that a certain antisymmetric
bilinear form is well--defined on the
Grothendieck group of a
cluster--tilted algebra
associated with a
finite--dimensional hereditary algebra.
The conjecture
was proved in~\cite{Pcc} in the more
general context of Hom-finite
2-Calabi--Yau triangulated
categories. It was used in order to
study the existence of a cluster
character on such a category $\cat$,
by using a formula proposed by
Caldero--Keller.

In the present paper,
we restrict to the case where $\cat$
is algebraic (i.e. is the stable
category of a Frobenius category).
We first use this bilinear form to
prove a generalized mutation rule
for quivers of cluster--tilting
subcategories in $\cat$.
When the cluster--tilting subcategories
are related by a single mutation,
this shows,
via the method of~\cite{GLSrigid},
that their quivers are related by
the Fomin--Zelevinsky mutation rule.
This special case was already proved
in~\cite{BIRS}, without
assuming $\cat$ to be
algebraic.

We also compute the
Grothendieck group of the
triangulated category $\cat$.
In particular, this
allows us to improve on
results by M. Barot, D. Kussin
and H. Lenzing: We compare the
Grothendieck group of a cluster
category $\cat_A$ with the group
$\kred$. The latter group was defined
in~\cite{BKL}
by only considering the triangles
in $\cat_A$ which are induced by
those of the derived category.
More precisely, we prove that
those two groups
are isomorphic for any cluster
category associated with a finite
dimensional hereditary algebra, with its
triangulated structure defined by
B. Keller in~\cite{Ktri}.

This paper is organized as follows:
The first section is dedicated to
notation and necessary background
from~\cite{FZ1}, \cite{GLSrigid},
\cite{KR1}, \cite{Pcc}.
In section~\ref{section: kzero},
we compute the Grothendieck group
of the triangulated category $\cat$.
In section~\ref{section: mutation},
we prove a generalized mutation rule
for quivers of cluster--tilting
subcategories in $\cat$. In particular,
this yields a new proof of the
Fomin--Zelevinsky mutation rule,
under the restriction that $\cat$
is algebraic. We finally show that
$\kzero(\cat_A)=\kred$ for any
finite dimensional hereditary
algebra $A$.

\section*{Acknowledgements}

This article is part of my PhD thesis, under the supervision of Professor
B. Keller. I would like to thank him deeply for introducing me to the subject 
and for his infinite patience. 

\tableofcontents

\section{Notations and background}

Let $\ec$ be a Frobenius category
whose idempotents split and which is
linear over a given algebraically
closed field $k$.
By a result of Happel~\cite{Happel}, its
stable category $\cat = \es$ is triangulated.
We assume moreover, that $\cat$ is Hom-finite,
2-Calabi--Yau and has a cluster--tilting
subcategory (see section~\ref{subsection: cts}),
and we denote by $\susp$ its
suspension functor. Note that we do not
assume that $\ec$ is Hom-finite.

We write $\xc(\;,\;)$, or
$\homph_{\xc}(\;,\;)$,
for the morphisms in
a category $\xc$ and $\homph_{\xc}(\;,\;)$
for the morphisms in the category of
$\xc$-modules. We also denote by
$X\hat\,$ the projective
$\xc$-module represented by $X$:
$X\hat\, = \xc(?,X)$.

\subsection{Fomin--Zelevinsky mutation for matrices}\label{subsection: FZ}

Let $B = (b_{ij})_{i,j\in I}$ be a
finite or infinite matrix, and let
$k$ be in $I$. The Fomin and Zelevinsky mutation
of $B$ (see~\cite{FZ1}) in direction $k$ is
the matrix
$$ \mu_k(B) = (b_{ij}')$$
defined by
$$
b_{ij}' = \left\{ \begin{array}{ll}
 -b_{ij} & \text{if } i=k \text{ or } j=k, \\
  b_{ij} + \frac{|b_{ik}|b_{kj} +b_{ik}|b_{kj}|}{2} & \text{else.}
\end{array}\right.
$$
Note that 
$\mu_k\big{(} \mu_k(B) \big{)} = B$ and that
if $B$ is skew-symmetric, then so is $\mu_k(B)$.

We recall two lemmas of~\cite{GLSrigid},
stated for infinite matrices,
which will be useful in section~\ref{section: mutation}.
Note that lemma 7.2 is a restatement of~\cite[(3.2)]{BFZ}.
Let $S=(s_{ij})$ be the matrix defined by
$$
s_{ij} = \left\{\begin{array}{ll}
 -\delta_{ij} + \frac{|b_{ij}| - b_{ij}}{2} & \text{if } i=k,\\
 \delta_{ij} & \text{else.}
\end{array}\right.
$$

\noindent{\bfseries Lemma 7.1~(\cite[Geiss--Leclerc--Schr\"oer]{GLSrigid})} :
{\itshape
Assume that $B$ is skew-symmetric. Then,
$S^2 = 1$ and the $(i,j)$-entry of the transpose
of the matrix $S$ is given by }
$$
s^\text{t}_{ij} = \left\{\begin{array}{ll}
-\delta_{ij} + \frac{|b_{ij}| + b_{ij}}{2} & \text{if } j=k,\\
 \delta_{ij} & \text{else.}
\end{array}\right.
$$

The matrix $S$ yields a convienent way
to describe the mutation of $B$
in the direction $k$:

\noindent{\bfseries Lemma 7.2~(\cite[Geiss--Leclerc--Schr\"oer]{GLSrigid},
~\cite[Berenstein--Fomin--Zelevinsky]{BFZ})} : 
{\itshape
Assume that $B$ is skew-symmetric. Then we have: }
$$
\mu_k(B) = S^\text{t} B S.
$$

Note that the product is well-defined
since the matrix $S$ has a finite number
of non vanishing entries in
each column.

\subsection{Cluster--tilting subcategories}\label{subsection: cts}

A cluster--tilting subcategory 
(see~\cite{KR1}) of $\cat$
is a full subcategory $\tc$ such that
\begin{itemize}
 \item[a)] $\tc$ is a linear subcategory;
 \item[b)] for any object $X$ in $\cat$, the contravariant 
functor $\cat( ? , X)|_{\tc}$ is finitely generated;
 \item[c)] for any object $X$ in $\cat$, 
we have $\cat(X,\susp T) = 0$ 
for all $T$ in $\tc$ if and only if $X$ belongs to $\tc$.
\end{itemize}

We now recall some results from~\cite{KR1}, 
which we will use in the sequel.
Let $\tc$ be a cluster--tilting subcategory of $\cat$, 
and denote by $\mc$ its preimage in $\ec$. 
In particular $\mc$ contains the full subcategory 
$\pc$ of $\ec$ formed by the projective-injective objects, 
and we have $\ms = \tc$.

The following proposition will be 
used implicitly, extensively in this paper.
\noindent {\bfseries Proposition~\cite[Keller--Reiten]{KR1} } : 
{\itshape
\begin{itemize}
 \item[a)] The category $\modms$ of finitely 
presented $\ms$-modules is abelian.
 \item[b)] For each object $X\in\cat$, there is a triangle 
$$
\susm X \gfl T^X_1 \gfl T^X_0 \gfl X
$$
of $\cat$, with $T^X_0$ and $T^X_1$ in $\tc$.
\end{itemize}
}

Recall that the perfect derived category $\perm$ is the 
full triangulated subcategory of the derived
category of ${\sf \mathcal{D}}\modm$ generated by 
the finitely generated projective $\mc$-modules.

\noindent {\bfseries Proposition~\cite[Keller--Reiten]{KR1} } : 
{\itshape
\begin{itemize}
 \item[a)] For each $X\in\ec$, there are conflations
$$
0 \gfl M_1 \gfl M_0 \gfl X \gfl 0 \;\text{ and }\;
0 \gfl X \gfl M^0 \gfl M^1 \gfl 0
$$
in $\ec$, with $M_0$, $M_1$, $M^0$ and $M^1$ in $\mc$.
 \item[b)] Let $Z$ be in $\modms$. Then $Z$ 
considered as an $\mc$-module lies in 
the perfect derived category $\perm$ 
and we have canonical isomorphisms
$$
D(\perm)(Z, ?) \simeq (\perm)(?,Z[3]).
$$
\end{itemize}
}

\subsection{The antisymmetric bilinear form}\label{subsection: fba}

In section~\ref{section: mutation}, 
we will use the existence of the 
antisymmetric bilinear form 
$\langle\;\,,\;\rangle_a$ on 
$\komodms$. We thus recall its definition 
from~\cite{CK1}. 

Let $\langle\;\,,\;\rangle$ be a truncated 
Euler form on $\modms$ defined by 
$$
\langle M,N \rangle = \dim \homph_{\ms}(M,N) - \dim \ext^1_{\ms}(M,N) 
$$
for any $M,N\in\modms$.
Define $\langle\;\,,\;\rangle_a$ to be the 
antisymmetrization of this form: 
$$
\langle M,N \rangle_a = \langle M,N \rangle - \langle N,M \rangle.
$$
This bilinear form descends to the 
Grothendieck group $\komodms$:

\noindent {\bfseries Lemma}~\cite[section 3]{Pcc} :
The antisymmetric bilinear form 
$$
\langle M,N \rangle_a : 
\komodms\times\komodms \gfl \zb
$$
is well-defined.

\section{Grothendieck groups of algebraic 2-CY categories 
with a cluster--tilting subcategory}\label{section: kzero}

We fix a cluster-tilting subcategory $\tc$ of $\cat$, 
and we denote by $\mc$ its preimage in $\ec$. 
In particular $\mc$ contains the full subcategory 
$\pc$ of $\ec$ formed by the projective-injective objects, 
and we have $\ms = \tc$.

We denote by $\hbe$ and $\dbe$ respectively the 
bounded homotopy category and the bounded derived category 
of $\ec$. We also denote by $\eace$, $\hbp$, $\hbm$ and $\eac$ 
the full subcategories of $\hbe$ whose objects are the 
$\ec$-acyclic complexes, the complexes of projective objects 
in $\ec$, the complexes of objects of $\mc$ 
and the $\ec$-acyclic complexes of objects of $\mc$, respectively.

\subsection{A short exact sequence of triangulated categories}\label{subsection: ses}

\begin{lem}\label{lem: trieq}
 Let $\ac_1$ and $\ac_2$ be thick, full triangulated subcategories of a triangulated
 category $\ac$ and let $\bc$ be $\ac_1 \cap \ac_2$.
Assume that for any object $X$ in $\ac$ there is a triangle
$X_1 \gfl X \gfl X_2 \gfl \susp X_1$ in $\ac$, with $X_1$ in $\ac_1$ 
and $X_2$ in $\ac_2$.
Then the induced functor $\ac_1 /\bc \gfl \ac/\ac_2 $
is a triangle equivalence.
\end{lem}

\begin{proof}
Under these assumptions, denote by $F$ the induced
triangle functor from $\ac_1 /\bc$ to $\ac/\ac_2$.
We are going to show that the functor $F$ is a full,
conservative, dense functor.
Since any full conservative triangle functor is
fully faithful, $F$ will then be an equivalence of categories.

We first show that $F$ is full.
Let $X_1$ and $X'_1$ be two objects in $\ac_1$.
Let $f$ be a morphism from $X_1$ to $X'_1$ in $\ac/\ac_2$
and let
$$
\xymatrix@-0.8pc{
 &  Y \bdr^{w} \bg & \\
X_1 & & X'_1
}
$$
be a left fraction which
represents $f$.
The morphism $w$ is
in the multiplicative system associated
with $\ac_2$
and thus
yields a triangle $\susm A_2 \fl Y \stackrel{w}{\gfl} X'_1 \fl A_2$
where $A_2$ lies in the subcategory $\ac_2$.
Moreover, by assumption, there exists
a triangle $Y_1 \fl Y \fl Y_2 \fl \susp Y_1$ with $Y_i$ in $\ac_i$.
Applying the octahedral axiom to the composition
$Y_1 \fl Y \fl X'_1$ yields a commutative diagram
whose two middle rows and columns are triangles in $\ac$
$$
\xymatrix{
 & \susm A_2 \bas \dreg & \susm A_2 \bas & \\
Y_1 \baseg \dr & Y \bas \dr & Y_2 \bas \dr & \susp Y_1 \baseg \\
Y_1 \dr & X'_1 \bas \dr & Z \bas \dr & \susp Y_1 \\
 & A_2 \dreg & A_2 & .
}
$$
Since $Y_2$ and $A_2$ belong to $\ac_2$, so does $Z$.
And since $X'_1$ and $Y_1$ belong to $\ac_1$, so does $Z$.
This implies, that
$Z$ belongs to $\bc$. The morphism $Y_1 \fl X'_1$ is
in the multiplicative system of
$\ac_1$ associated with $\bc$
and the diagram
$$
\xymatrix@-0.8pc{
 &  Y_1 \bdr \bg & \\
X_1 & & X'_1
}
$$
is a left fraction which represents $f$. This implies that $f$
is the image of a morphism in $\ac_1 / \bc$. Therefore
the functor $F$ is full.

We now show that $F$ is conservative.
Let $X_1 \stackrel{f}{\gfl} Y_1 \fl Z_1 \fl \susp X_1$ be
a triangle in $\ac_1$. Assume that $Ff$ is
an isomorphism in $\ac /  \ac_2$, which implies that
$Z_1$ is an
object of $\ac_2$. Therefore,
$Z_1$ is an
object of $\bc$ and $f$ is an
isomorphism in $\ac_1 / \bc$.

We finally show that $F$ is dense.
Let $X$ be an object of
the category $\ac / \ac_2$,
and let $X_1 \fl X \fl X_2 \fl \susp X_1$
be a triangle in $\ac$ with $X_i$
in $\ac_i$.
Since $X_2$ belongs to $\ac_2$, the image
of the morphism $X_1 \fl X$ in $\ac / \ac_2$
is an isomorphism. Thus $X$ is isomorphic to
the image by $F$ of an object in $\ac_1 / \bc$.
\end{proof}

As a corollary, we have the following:
\begin{lem}\label{lem:ses}
The following sequence of triangulated categories is short exact:
$$
0 \gfl \eac \gfl \hbm \gfl \dbe \gfl 0.
$$
\end{lem}

\noindent Remark: 
This lemma remains true if $\cat$
 is $d$-Calabi--Yau and $\ms$ is $(d-1)$-cluster--tilting, using 
 section 5.4 of~\cite{KR1}.

\begin{proof}
For any object $X$ in $\hbe$,
the existence of an object $M$
in $\hbm$ and of a quasi-isomorphism $w$
from $M$ to $X$ is obtained using the
approximation conflations of
Keller--Reiten (see section~\ref{subsection: cts}).
Since the cone of the morphism $w$
belongs to $\eace$, lemma~\ref{lem: trieq}
applies to the subcategories
$\eac$, $\hbm$ and $\eace$ of
$\hbe$.
\end{proof}

\begin{prop}\label{prop: diag}
$\!\!$The following diagram is commutative with exact rows and columns:
 $$\xymatrix@-3pt{
 & & 0 & 0 & & \\
0 \dr & \eac \dr^{i_{\mc}\;\;\;\;\;} \ar@{=}[d] & \hbm/\hbp \ha \dr & \es \dr \ha  & 0 & \\
0 \dr & \eac \dr & \hbm \dr \ha & \dbe \dr \ha & 0 & (D)  \\
 & & \hbp  \ar@{=}[r]  \ha & \hbp \ha_{i_{\pc}} \dr & 0 & \\
 & & 0 \ha & 0 \ha & & .
}$$
\end{prop}

\begin{proof}
 The column on the right side has been shown to be exact in~\cite{KV} 
and~\cite{Ri91}. The second row is exact by lemma~\ref{lem:ses}. 
The subcategories $\eac$ and $\hbp$ of $\hbm$ are left and right orthogonal 
to each other. This implies that the induced functors $i_{\mc}$ 
and $i_{\pc}$ are fully faithful and that 
taking the quotient of $\hbm$ by those two subcategories either in one order 
or in the other gives the same category. Therefore the first row is 
exact.
\end{proof}

\subsection{Invariance under mutation}\label{subsection: invariance}

A natural question is then to which extent the diagram $(D)$ depends on the 
choice of a particular cluster--tilting subcategory. Let thus 
$\tc'$ be another cluster--tilting subcategory of $\cat$, and let
$\mc'$ be its preimage in $\ec$. 
Let $\modm$ (resp. $\modm'$) be the category of $\mc$-modules 
(resp. $\mc'$-modules), i.e. of $k$-linear contravariant functors 
from $\mc$ (resp. $\mc'$) to the category of 
$k$-vector spaces. 

Let $X$ be the $\mc$-$\mc'$-bimodule which sends the 
pair of objects $(M,M')$ to the $k$-vector space 
$\ec(M,M')$. The bimodule $X$ induces a functor 
$F = ?\otimes_{\mc'}X : \modm' \gfl \modm$ 
denoted by $T_X$ in~\cite[section 6.1]{Kddc}.

Recall that the perfect derived category $\perm$ is the 
full triangulated subcategory of the derived
category ${\sf \mathcal{D}}\modm$ generated by 
the finitely generated projective $\mc$-modules.

\begin{prop}\label{prop:eq}
 The left derived functor 
$$
\begin{array}{crcl}
\lb F : & \der' & \gfl & \der
\end{array}
$$
is an equivalence of categories.
\end{prop}

\begin{proof}
Recall that if $X$ is an object in a category $\xc$, we denote by 
$X\hat\,$ the functor $\xc(?,X)$ represented by $X$.
 By \cite[6.1]{Kddc}, it is enough to check the 
following three properties:
\begin{itemize}
 \item[1.] For all objects $M'$, $M''$ of $\mc$, the group 
$\homph_{\der} \left( \lf M'\hat \, ,\, \lf  M''\hat\,  [n] \right)$ 
vanishes for $n \neq 0$ and identifies with 
$\homph_{\mc'} \left(M', M'' \right)$ for $n=0$;
 \item[2.] for any object $M'$ of $\mc'$, the complex $\lf M'\hat\,$ belongs to $\perm$;
 \item[3.] the set $\left\{ \lf M'\hat\,  , \,  M' \in \mc' \right\}$ generates $\der$ 
as a triangulated category with infinite sums.
\end{itemize}

Let $M'$ be an object of $\mc'$, and let 
$\xymatrix@-6pt{M_1 \; \dri & M_0 \dre & M'}$ be a conflation in $\ec$, 
with $M_0$ and $M_1$ in $\mc$, and whose deflation is a 
right $\mc$-approximation (c.f. section~4 of~\cite{KR1}). 
The surjectivity of the map
$M_0\hat\, \gfl  \ec(?,M')|_{\mc}$ implies that the complex 
$P = (\cdots \fl 0 \fl M_1\hat\, \fl M_0\hat\, \fl 0 \fl \cdots)$ 
is quasi-isomorphic to $\lf M'\hat\, = \ec(?,M')|_{\mc}$.
Therefore $\lf M'\hat\,$ belongs to the subcategory 
$\perm$ of $\der$. Moreover, we have, for any $n\in\zb$ 
and any $M''\in\mc'$, the equality
$$
\homph_{\der} \left( \lf M'\hat\,,\lf M''\hat\,[n] 
\right) = \homph_{\hc^\text{b} \modm}\left(
P, \ec(?,M'')|_{\mc}[n] \right)
$$
where the right--hand side vanishes for $n\neq 0,1$.
In case $n=1$ it also vanishes, since 
$\ext^1_{\ec}(M',M'')$ vanishes. Now, 
\begin{eqnarray*}
\homph_{\hc^\text{b} \modm}\left(
P, \ec(?,M'')|_{\mc} \right) & \simeq & 
\ke \left( \ec(M_0,M'') \fl \ec(M_1,M'') \right) \\
 & \simeq & \ec(M',M'').
\end{eqnarray*}

It only remains to be shown that the set 
$R = \left\{ \lf M'\hat\,  , \,  M' \in \mc' \right\}$ 
generates $\der$. Denote by $\rc$ the full triangulated 
subcategory with infinite sums of $\der$ generated 
by the set $R$. 
The set $\left\{M\hat\,, \,M\in\mc \right\}$
generates $\der$ as a triangulated category 
with infinite sums. Thus it is enough to show 
that, for any object $M$ of $\mc$, the complex 
$M\hat\,$ concentrated in degree $0$ belongs 
to the subcategory $\rc$.
Let $M$ be an object of $\mc$, and let 
$\xymatrix@-6pt{M\; \dri^{i} & M'_0 \dre^{p} & M'_1}$ 
be a conflation of $\ec$ with $M'_0$ and
$M'_1$ in $\mc'$. 
Since $\ext_{\ec}^1 (?,M)|_{\mc}$ vanishes, 
we have a short exact sequence of 
$\mc$-modules
$$
0\gfl \ec(?,M)|_{\mc} \gfl \ec(?,M'_0)|_{\mc} \gfl\ec(?,M'_1)|_{\mc} \gfl 0,
$$
which yields the triangle 
$$
M\hat\, \gfl \lf M'_{0}\hat\,\gfl \lf M'_{1}\hat\,\gfl \susp M\hat\,.
$$
\end{proof}

As a corollary of proposition~\ref{prop:eq}, up to equivalence 
the diagram $(D)$ does not depend on the choice of a 
cluster--tilting subcategory. To be more precise: 
Let $G$ be the functor which sends an object $X$ 
in the category $\hbmp$ to a representative of 
$(\lf) X\hat\,$ in $\hbm$, 
and a morphism in $\hbmp$ to the induced one in $\hbm$.

\begin{cor}\label{cor: diag}
 The following diagram is commutative
$$
\xymatrix{
& \derp \ddr^{\lf} & & \der \\
\hbmp \bbas \hdrm \ddr^{\;\;\;\;\;\;\;\;\;\;\;\;\;\;\;\;\;\;\;\;\;\;\;\;\; G}
& \ha & 
\hbm \bbas \hdrm & \\
 & \hbp \hagm \ar@{^{(}-}[u] \dreg \bgm & \dreg & \hbp \hagm \hham \bgm \\
\dbe \ddreg & & \dbe &
}
$$
and the functor $G$ is an equivalence of categories.
\end{cor}

We denote by $\perms$ the full subcategory 
of $\perm$ whose objects are the complexes 
with homologies in $\modms$. 
The following lemma will allow 
us to compute the Grothendieck group 
of $\perms$ in section~\ref{subsection: ko}:
\begin{lem}\label{lem: t-str}
The canonical t-structure on 
$\der$ restricts to a t-structure on 
$\perms$, whose heart is $\modms$.
\end{lem}

\begin{proof}
 
By~\cite{KVaisles}, it is enough to 
show that for any object 
$M^\bullet$ of $\perms$, 
its truncation $\tau_{\leq 0}M^\bullet$ 
in $\der$ belongs to $\perms$. 
Since $M^\bullet$ is in $\perms$, 
$\tau_{\leq 0}M^\bullet$ is bounded, 
and is thus formed from the complexes 
H$^i(M^\bullet)$ 
concentrated in one degree 
by taking iterated extensions. 
But, for any $i$, the $\mc$-module 
H$^i(M^\bullet)$ actually is an 
$\ms$-module. Therefore, 
by~\cite{KR1} 
(see section~\ref{subsection: cts}), 
it is perfect as an $\mc$-module 
and it lies in $\perms$.
\end{proof}

The next lemma already appears in~\cite{Tab}.
For the convenience of the reader, 
we include a proof.
\begin{lem}\label{lem: pereac}
 The Yoneda equivalence of triangulated categories $\hbm \gfl \perm$ 
induces a triangle equivalence $\eac \gfl \perms$.
\end{lem}

\begin{proof}
 We first show that the cohomology groups of an $\ec$-acyclic bounded complex $M$ 
vanish on $\pc$. Let $P$ be a projective object in $\ec$ and let $E$ 
be a kernel in $\ec$ of the map $M^n \gfl M^{n+1}$. Since $M$ is $\ec$-acyclic, 
such an object exists, and moreover, it is an image of the map $M^{n-1} \gfl M^n$.
Any map from $P$ to $M^n$ whose composition with $M^n \fl M^{n+1}$ vanishes 
factors through the kernel $E \infl M^n$. Since $P$ is projective, this factorization 
factors through the deflation $M^{n-1} \defl E$. 
$$
\xymatrix{
 & & P \bggp \ar@/_1.6pc/@{-->}[ddl] \bas \bddr^0 & & \\
M^{n-1} \bdre \ar@{-}[r] & \dr & M^n \ddr & & M^{n+1} \\
 & E \hdri & & &
}
$$
Therefore, we have $\coh^n(M\hat\,)(P) = 0$ for all projective objects $P$, and 
$\coh^n(M\hat\,)$ belongs to $\modms$. Thus the Yoneda functor 
induces a fully faithful functor from $\eac$ to $\perms$.
To prove that it is dense, it is enough 
to prove that any object of the heart $\modms$ of the t-structure 
on $\perms$ is in its essential image.

But this was proved in~\cite[section 4]{KR1} 
(see section~\ref{subsection: cts}).

\end{proof}

\begin{prop}\label{prop: perms}
 There is a triangle equivalence of categories 
$$
\perms \stackrel{\simeq}{\gfl} \permsp
$$
\end{prop}

\begin{proof}
Since the categories $\hbp$ and
$\eac$ are left-right orthogonal
in $\hbm$, this is immediate from
corollary~\ref{cor: diag}
and lemma~\ref{lem: pereac}.
\end{proof}

\subsection{Grothendieck groups}
\label{subsection: ko}

For a triangulated (resp. additive, resp. abelian) category $\ac$, 
we denote by $\kzero^\text{tri}(\ac)$ or simply $\kzero(\ac)$ 
(resp. $\kzero^\text{add}(\ac)$, resp. $\kzero^\text{ab}(\ac)$) 
its Grothendieck group (with respect 
to the mentioned structure of the category). 
For an object $A$ in $\ac$, we also denote by $[A]$ 
its class in the Grothendieck group of $\ac$. 

The short exact sequence of triangulated categories
$$
0 \gfl \eac \gfl \hbm/\hbp \gfl \es \gfl 0
$$
given by proposition~\ref{prop: diag} 
induces an exact sequence in the Grothendieck groups 
$$
(\ast) \;\;\;\;\; \kzero\big{(} \eac \big{)} \gfl \kzero\big{(} 
\hbm/\hbp \big{)} \gfl \kzero\big{(} \es \big{)} \gfl 0 .
$$

\begin{lem}\label{lem: exseq}
 The exact sequence $(\ast)$ is isomorphic to an exact sequence
$$
(\ast\ast) \;\;\;\;\; \kzero^\text{ab}\big{(} \modms 
\big{)} \stackrel{\varphi}{\gfl} \kzero^\text{add}\big{(} \ms \big{)} 
\gfl \kzero^\text{tri}\big{(} \es \big{)} \gfl 0 .
$$
\end{lem}

\begin{proof}
 
First, note that, by~\cite{Tab}, 
see also lemma~\ref{lem: pereac}, we have an isomorphism 
between the Grothendieck groups
$\kzero\big{(} \eac \big{)}$ and $\kzero\big{(} \perms \big{)}$.
The t-structure on $\perms$ whose heart is $\modms$, 
see lemma~\ref{lem: t-str}, in turn yields 
an isomorphism between the
Grothendieck groups $\kzero^\text{tri}\big{(} \perms \big{)}$ 
and $\kzero^\text{ab}\big{(} \modms \big{)}$.
Next, we show that the canonical additive functor 
$\ms \stackrel{\al}{\gfl} \hbm/\hbp$ induces an isomorphism 
between the Grothendieck groups 
$\kzero^\text{add}\big{(} \ms \big{)}$ and
$\kzero^\text{tri}\big{(} \hbm/\hbp \big{)}$. 
For this, let us consider the canonical additive 
functor $\ms \stackrel{\beta}{\gfl} \hbms$ and 
the triangle functor 
$\hbm \stackrel{\gamma}{\gfl}\hbms$.
The following diagram describes 
the situation: 
$$
\xymatrix{
\hbms & \hbm \ar@{->}[l]_{\gamma} \bas \\
\ms \ha^{\beta} \dr^{\al\;\;\;\;\;\;\;\;\;} & 
\hbm/\hbp \hagp_{\gamma}
}
$$
The functor $\gamma$ vanishes on the full 
subcategory $\hbp$, thus 
inducing a triangle functor, still denoted by $\gamma$, 
from $\hbm/\hbp$ to $\hbms$.
Furthermore, the functor $\beta$ induces an isomorphism at the level of 
Grothendieck groups, whose inverse $\kzero(\beta)^{-1}$ is given by 
\begin{eqnarray*}
\kzero^\text{tri} \big{(} \hbms \big{)} & \gfl & \kzero^\text{add} \big{(} \ms \big{)} \\ 
\, [M] & \longmapsto & \sum_{i\in\zb} (-1)^i [M^i].
\end{eqnarray*}
As the group $\kzero^\text{tri}\big{(} \hbm/\hbp \big{)}$ is generated by 
objects concentrated in degree $0$, it is straightforward to 
check that the morphisms $\kzero(\al)$ and 
$\kzero (\beta)^{-1} \kzero (\gamma)$ are inverse 
to each other.

\end{proof}

As a consequence of the exact sequence $(\ast\ast)$,
we have an isomorphism between
$\kzero^\text{tri}\big{(} \es \big{)}$ and
$\kzero^\text{add}\big{(} \ms \big{)}/\im \varphi$.
In order to compute $\kzero^\text{tri}\big{(} \es \big{)}$,
the map $\varphi$ has to be made explicit.
We first recall some results from Iyama--Yoshino~\cite{IY}
which generalize results from~\cite{BMRRT}:
For any indecomposable $M$ of $\ms$ not in $\pc$,
there exists $M^\ast$ unique up to isomorphism
such that $(M,M^\ast)$ is an exchange pair.
This means that $M$ and $M^\ast$ are not isomorphic and that
the full additive subcategory of
$\cat$ generated by all the indecomposable objects
of $\ms$ but those isomorphic to $M$, and by the
indecomposable objects isomorphic to $M^\ast$ is
again a cluster--tilting subcategory.
Moreover, $\dim \es(M,\susp M^\ast)=1$.
We can thus fix two (non-split)
exchange triangles
$$
M^\ast \fl B_M \fl M \fl \susp M^\ast \text{ and }
M \fl B_{M^\ast} \fl M^\ast \fl \susp M.
$$ 
We may now state the following:
\begin{theo}\label{theo: kzero}
 The Grothendieck group of the triangulated category $\es$ is 
the quotient of that of the additive subcategory $\ms$ by all 
relations $[B_{M^\ast}] - [B_M]$: 
$$
\kzero^\text{tri}\big{(} \es \big{)} \simeq 
\kzero^\text{add}\big{(} \ms \big{)}/([B_{M^\ast}] - [B_M])_M.
$$
\end{theo}

\begin{proof}
We denote by $S_M$ the simple $\ms$-module associated to 
the indecomposable object $M$. 
This means that $S_M(M')$ vanishes for all indecomposable
objects $M'$ in $\ms$ not isomorphic to $M$ and that
$S_M(M)$ is isomorphic to $k$. The abelian group
$\kzero^\text{ab}\big{(}\modms\big{)}$ is generated
by all classes $[S_M]$. In view of lemma~\ref{lem: exseq},
it is sufficient to prove that the image of the class
$[S_M]$ under $\varphi$ is $[B_{M^\ast}] - [B_M]$.
First note that the $\mc$-module $\ext^1_{\ec}(?,M^\ast)|_{\mc}$
vanishes on the projectives ; it can thus be viewed as an
$\ms$-module, and as such, is isomorphic to $S_M$.
After replacing $B_M$ and $B_{M'}$ by isomorphic objects of $\es$,
we can assume that the exchange triangles
$M^\ast \fl B_M \fl M \fl \susp M^\ast$ and
$M \fl B_{M^\ast} \fl M^\ast \fl \susp M$
come from conflations
$\xymatrix@-6pt{M^\ast\;\, \dri & B_M \dre & M}$ and
$\xymatrix@-6pt{M\;\, \dri & B_{M^\ast} \dre & M^\ast.}$
The spliced complex
$$(\cdots \fl 0 \fl M \fl B_{M^\ast}
\fl B_M \fl M \fl 0 \fl \cdots)$$ denoted by $C^\bullet$,
 is then an $\ec$-acyclic complex,
and it is the image of $S_M$ under the functor
$\modms \subset \perms \simeq \eac$. Indeed, we have two
long exact sequences induced by the conflations above:
$$0 \fl \mc(?,M) \fl \mc(?,B_{M^\ast})
\fl \ec(?,M^\ast)|_{\mc} \fl \ext^1_{\ec}(?,M)|_{\mc} = 0 \text{ and}
$$
$$
0 \fl \ec(?,M^\ast)|_{\mc} \fl \mc(?,B_M) \fl \mc(?,M)
\fl \ext^1_{\ec}(?,M^\ast)|_{\mc} \fl \ext^1_{\ec}(?,B_M)|_{\mc}.$$
Since $B_M$ belongs to $\mc$, the functor
$\ext^1_{\ec}(?,B_M)$ vanishes on $\mc$, and the complex:
$$(C\hat\,) :\;\;\;\;\;\; (\cdots \fl 0 \fl M\hat\, \fl (B_{M^\ast})\hat\,
\fl (B_M)\hat\, \fl M\hat\, \fl 0 \fl \cdots)$$
is quasi-isomorphic to $S_M$.

Now, in the notations of the proof of lemma~\ref{lem: exseq}, 
$\varphi[S_M]$ is the image of the class of the 
$\ec$-acyclic complex complex $C^\bullet$ under 
the morphism $\kzero (\beta)^{-1} \kzero (\gamma)$. 
This is $[M]- [B_M] + [B_{M^\ast}] - [M]$ 
which equals $[B_{M^\ast}]-[B_M]$ as claimed.
\end{proof}

\section{The generalized mutation rule}\label{section: mutation}

Let $\tc$ and $\tc'$ be two
cluster--tilting subcategories
of $\cat$. Let $Q$ and $Q'$ be
the quivers obtained from their
Auslander--Reiten quivers
by removing all loops and
oriented 2-cycles.

Our aim, in this section, is to
give a rule relating $Q'$ to $Q$,
and to prove that it generalizes
the Fomin--Zelevinsky mutation
rule.

\noindent {\itshape Remark: }
\begin{itemize}
 \item[.] Assume that $\cat$ has cluster--tilting objects.
Then it is proved in~\cite[Theorem I.1.6]{BIRS}, without
assuming that $\cat$ is algebraic,
that the Auslander--Reiten quivers
of two cluster--tilting objects having
all but one indecomposable direct summands
in common (up to isomorphism) are related
by the Fomin--Zelevinsky mutation
rule.
 \item[.] To prove that the generalized
mutation rule actually generalizes
the Fomin--Zelevinsky mutation
rule, we use the ideas of
section 7 of~\cite{GLSrigid}.
\end{itemize}

\subsection{The rule}
\label{subsection: gmr}

As in section~\ref{section: kzero}, we fix
a cluster--tilting subcategory $\tc$ of $\cat$,
and write $\mc$ for its preimage in $\ec$,  so that
$\tc = \ms$.
Define $Q$ to be the quiver obtained from
the Auslander--Reiten quiver
of $\ms$ by deleting its loops
and its oriented 2-cycles.
Its vertex corresponding to an indecomposable
object $L$ will also be labeled by $L$.
We denote by $a_{LN}$ the number of arrows
from vertex $L$ to vertex $N$ in the quiver $Q$.
Let $B_{\mc}$ be the matrix whose entries are
given by $b_{LN} = a_{LN} - a_{NL}$.

Let $R_{\mc}$ be the matrix of
$\langle \; , \; \rangle_a : \komodms\times\komodms \gfl \zb$
in the basis given by the classes of the simple
modules.

\begin{lem}\label{lem: rmegalbm}
 The matrices $R_{\mc}$ and $B_{\mc}$ are equal: 
$R_{\mc} = B_{\mc}$.
\end{lem}

\begin{proof}
 Let $L$ and $N$ be two non-projective indecomposable 
objects in $\mc$. Then
$\dim \homph (S_L, S_N) - \dim \homph (S_N, S_L) = 0$ and 
we have: 
$$\langle [S_L], [S_N] \rangle_a = \dim \ext^1(S_N,S_L) - 
\dim\ext^1(S_L,S_N) = b_{L,N}.$$
\end{proof}

Let $\tc'$ be another cluster--tilting subcategory
of $\cat$, and let $\mc'$ be its preimage in the 
Frobenius category $\ec$.
Let $(M'_i)_{i\in I}$
(resp. $(M_j)_{j\in J}$)
be representatives
for the isoclasses of non-projective
indecomposable objects in $\mc'$
(resp. $\mc$).
The equivalence of categories
$\perms \stackrel{_\sim}{\gfl} \permsp$
 of proposition~\ref{prop: perms}
induces an isomorphism
between the Grothendieck groups $\komodms$
and $\komodmsp$ whose matrix,
in the bases given by the classes
of the simple modules,
is denoted by $S$.
The equivalence of categories
$\der \stackrel{_\sim}{\gfl} \derp$ restricts
to the identity on $\hbp$,
so that it induces an equivalence
$\perm/\perP \stackrel{_\sim}{\gfl} \permp/\perP$.
Let $T$ be the matrix of the induced
isomorphism from $\koprojmsurp$
to $\koprojmpsurp$, in the bases
given by the classes
$[\mc(?,M_j)]$, $j\in J$,
and $[\mc'(?,M'_i)]$, $i\in I$.
The matrix $T$ is much easier
to compute than the matrix $S$.
Its entries $t_{i j}$ are given
by the approximation triangles of
Keller and Reiten in the following
way: For all $j$, there exists a
triangle of the form
$$
\susm M_j \gfl
\bigoplus_i \beta_{ij} M'_i \gfl
\bigoplus_i \al_{ij} M'_i \gfl M_j.
$$
Then, we have:

\begin{theo}
\label{theo: gmr}
 \begin{itemize}
 \item[a)] \emph{(Generalized mutation rule)} 
The following equalities hold:
$$t_{i j} = \al_{ij} - \beta_{ij}$$
and 
$$B_{\mc'} = T B_{\mc} T^\text{\emph{t}}.$$
 \item[b)] The category $\cat$ has a 
cluster--tilting object if and only if 
all its cluster--tilting subcategories 
have a finite number of pairwise non-isomorphic 
indecomposable objects.
 \item[c)] All cluster--tilting objects 
of $\cat$ have the same number of indecomposable 
direct summands (up to isomorphism).
\end{itemize}
\end{theo}

Note that point c) was shown in~\cite[5.3.3(1)]{IAuscorr} 
(see also~\cite[I.1.8]{BIRS}) and, in a more general 
context, in~\cite{DK}.
Note also that, for the generalized mutation rule
to hold, the cluster--tilting subcategories
do not need to be related by a sequence of mutation.

\begin{proof}
 Assertions b) and c) are consequences of
the existence of an isomorphism between
the Grothendieck groups $\komodms$
and $\komodmsp$. Let us prove
the equalities a).
Recall from~\cite[section 3.3]{Pcc},
that the antisymmetric bilinear form
$\langle \;\, ,\;\rangle_a$ on $\modms$
is induced by the usual Euler form
$\langle \;\, ,\;\rangle_E$ on
$\perms$. The following commutative
diagram
$$
\xymatrix{
\perms\times\perms \bdr_{\langle \;,\;\rangle_E}
\ddr^{\simeq} & & \permsp\times\permsp
\bg^{\langle \;,\;\rangle_E} \\
& \zb & ,
}
$$
thus induces a commutative diagram
$$
\xymatrix{
\komodms\times\komodms \bdr_{\langle \;,\;\rangle_a}
\ddr^{S\times S} & & \komodmsp\times\komodmsp
\bg^{\langle \;,\;\rangle_a} \\
& \zb & .
}
$$
This proves the equality
$R_{\mc} = S^\text{t} R_{\mc'}S$,
or, by lemma~\ref{lem: rmegalbm},
$$
\!\!\!\!\!\!\!\!\!\!\!\!\!\!\!\!\!\!\!\!\!
(1) \;\;\;\;\;\;\; B_{\mc} = S^\text{t} B_{\mc'}S.
$$
Any object of $\perms$ becomes
an object of $\perm/\perP$ through
the composition
$
\perms \hookrightarrow \perm \defl
\perm/\perP.
$
Let $M$ and $N$ be two non-projective
indecomposable objects in $\mc$.
Since $S_N$ vanishes on $\pc$,
we have
\begin{eqnarray*}
 \homph_{\perm/\perP}\big(
\mc(?,M),S_N\big) & = &
\homph_{\perm}\big(
\mc(?,M),S_N\big) \\
 & = & \homph_{\modm}\big(
\mc(?,M),S_N\big)\\
 & = & S_N(M).
\end{eqnarray*}
Thus $\dim\homph_{\perm/\perP}\big(
\mc(?,M),S_N\big) = \delta_{M N}$,
and the commutative diagram
$$
\xymatrix{
\perm/\perP\times\perm/\perP \bdr_{R\Hom}
\ddr^{\simeq} & & \permp/\perP\times\permp/\perP
\bg^{R\Hom} \\
& \operatorname{per} k & ,
}
$$
induces a commutative diagram
$$
\xymatrix{
\scriptstyle{ \koprojmsurp\times\komodms} \bdr_{\text{Id}}
\ddr^{T\times S} & & {\scriptstyle \koprojmpsurp\times\komodmsp}
\bg^{\text{Id}} \\
& \zb & .
}
$$
In other words, the matrix $S$ is
the inverse of the transpose of $T$:
$$
\!\!\!\!\!\!\!\!\!\!\!\!\!\!\!\!\!\!
(2) \;\;\;\;\;\; S = T^{\text{-t}}
$$
Equalities (1) and (2) imply what
was claimed, that is
$$B_{\mc'} = T B_{\mc} T^\text{\emph{t}}.$$

Let us compute the matrix $T$: 
Let $M$ be indecomposable non-projective
in $\mc$, and let
$$
\susm M \gfl M_1' \gfl M_0' \gfl M
$$
be a Keller--Reiten approximation triangle
of $M$ with respect to $\mc'$, which we may
assume to come from a conflation in $\ec$.
This conflation yields a projective
resolution
$$
0\gfl (M_1')\hat\,\gfl (M_0')\hat\,\gfl
\ec(?,M)|_{\mc'} \gfl \ext^1_{\ec}(?,M_1')|_{\mc'} = 0.
$$
so that $T$ sends the class of
$M\hat\,$ to $[(M_0')\hat\,] - [(M_1')\hat\,]$.
Therefore, $t_{ij}$ equals $\al_{ij} - \beta_{ij}$.
\end{proof}

\newpage

\subsection{Examples}

\subsubsection{}

As a first example, let
$\cat$ be the cluster category
associated with the quiver
of type $A_4$: $1\fl2\fl3\fl4$.
Its Auslander--Reiten quiver
is the Moebius strip:
$$
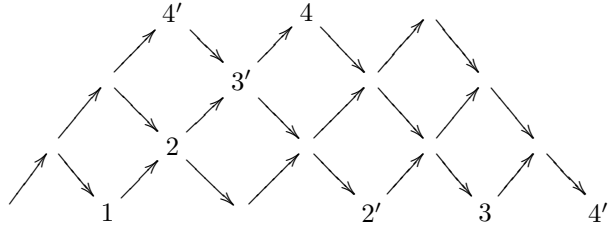
\xymatrix@-1pc{
 & & & 
4' \bdr & & 4 \bdr & & \bdr & & & \\
 & &  \bdr \hdr & & 3' \bdr \hdr & &  \bdr \hdr & & \bdr & & \\
 &  \bdr \hdr & & 2 \bdr \hdr & &  \bdr \hdr & &  \bdr \hdr & &  \bdr & \\
 \hdr & & 1 \hdr & &  \hdr & & 2' \hdr & 
& 3 \hdr & & 4'.
}
$$
Let $M=M_1\oplus M_2\oplus
M_3\oplus M_4$, where the indecomposable
$M_i$ corresponds to the vertex
labelled by $i$ in the picture.
Let also $M'=M'_1\oplus M'_2\oplus
M'_3\oplus M'_4$, where $M'_1=M_1$,
and where the indecomposable 
$M'_i$ corresponds to the vertex 
labelled by $i'$ if $i\neq1$.
One easily computes the following
Keller--Reiten approximation 
triangles:\\
$\susm M_1\gfl \;\;0\;\;\!\gfl M'_1\gfl M_1$,\\ 
$\susm M_2\gfl M'_2\gfl M'_1\gfl M_2$,\\
$\susm M_3\gfl M'_4\gfl \;\;0\;\;\!\gfl M_4$ and\\
$\susm M_4\gfl M'_4\gfl M'_3\gfl M_4$;\\
so that the matrix $T$ is given by:
$$
T = \left( \begin{array}{cccc}
1 & 1 & 0 & 0 \\
0 & -1\;\;\, & 0 & 0 \\
0 & 0 & 0 & 1 \\
0 & 0 & -1\;\;\, & -1\;\;\,
\end{array}\right).
$$
We also have
$$
B_{M'} = \left( \begin{array}{cccc}
0 & -1\;\;\, & 1 & 0 \\
1 & 0 & -1\;\;\, & 0 \\
-1\;\;\, & 1 & 0 & -1\;\;\, \\
0 & 0 & 1 & 0
\end{array}\right).
$$
Let maple compute 
$$
T^\text{\emph{-1}} B_{M'} T^\text{\emph{-t}} = 
\left( \begin{array}{cccc}
0 & 1 & 0 & 0 \\
-1\;\;\, & 0 & -1\;\;\, & 1 \\
0 & 1 & 0 & -1\;\;\, \\
0 & -1\;\;\, & 1 & 0
       \end{array}\right),
$$
which is $B_M$.

\subsubsection{}

Let us look at a more interesting
example, where one cannot easily
read the quiver of $M'$ from the
Auslander--Reiten quiver of $\cat$.
Let $\cat$ be the cluster category
associated with the quiver $Q$:
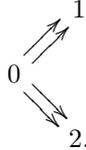
$$\xymatrix@-1pc{
& 1 \\
0 \hdr<.5ex> \hdr<-.5ex> 
 \bdr<.5ex> \bdr<-.5ex> & \\
& 2.
}$$
For $i=0,1,2$, let $M_i$ be (the image
in $\cat$ of) the projective indecomposable
(right) $kQ$-module associated with
vertex $i$. Their dimension vectors
are respectively $[1,0,0], [2,1,0]$
and $[2,0,1]$. Let $M$ be the
direct sum $M_0\oplus M_1\oplus M_2$.
Let $M'$ be the direct sum
$M'_0\oplus M'_1\oplus M'_2$,
where $M'_0, M'_1$ and $M'_2$
are (the images in $\cat$
of) the indecomposable regular
$kQ$-modules with dimension
vectors $[1,2,0], [0,1,0]$
and $[2,4,1]$ respectively.
As one can check, using~\cite{Kjava},
$M$ and $M'$ are two
cluster--tilting objects
of $\cat$. To compute
Keller--Reiten's approximation 
triangles, amounts to computing
projective resolutions in
$\!\!\!\!\mod kQ$, viewed as
$\!\!\!\!\mod\!\! \operatorname{End}_{\cat}(M)$.
One easily computes these projective
resolutions, by considering
dimension vectors:\\
$0 \gfl 8M_0\gfl M_2\oplus4M_1\gfl M'_2\gfl 0$,\\
$0 \gfl 2M_0\gfl M_1\gfl M'_1\gfl 0$ and\\
$0 \gfl 3M_0\gfl 2M_1\gfl M'_0\gfl 0$.\\
By applying the generalized mutation rule,
one gets the following quiver
$$
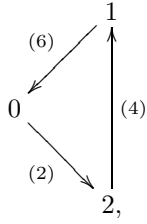
\xymatrix{
& 1 \bg_{(6)} \\
0 \bdr_{(2)} & \\
& 2 \ar@{->}[uu]_{(4)},
}$$
which is therefore the quiver of
$\operatorname{End}_{\cat}(M')$
since by~\cite{BMRclustermutation},
there are no loops or $2$-cycles
in the quiver of the endomorphism algebra
of a cluster--tilting object in
a cluster category.

\subsection{Back to the mutation rule}

We assume in this section that the
Auslander--Reiten quiver
of $\tc$ has no loops nor 2-cycles.
Under the notations of section~\ref{subsection: gmr},
let $k$ be in $I$ and let
 $(M_k,M_k')$ be an exchange pair
(see section~\ref{subsection: ko}).
We choose $\ms'$ to be the cluster-tilting 
subcategory of $\cat$ obtained from $\ms$ 
by replacing $M_k$ by $M_k'$, so that
$M_i' = M_i$ for all $i\neq k$. 
Recall
that $T$ is the matrix of the isomorphism
$\koprojmsurp \gfl \koprojmpsurp$.

\begin{lem}\label{lem: S}
Then, the $(i,j)$-entry of the
matrix $T$ is given by 
$$
t_{i j} = \left\{ \begin{array}{ll}
 -\delta_{i j} + \frac{|b_{i j}| + b_{i j}}{2} & \text{if } j=k \\
 \;\;\delta_{i j} & \text{else.}
\end{array} \right.
$$
\end{lem}

\begin{proof}
Let us apply theorem~\ref{theo: gmr} to
compute the matrix $T$.
For all $j\neq k$, the triangle
$ \susm M_j \fl 0 \fl M_j' = M_j$
is a
Keller--Reiten approximation triangle
of $M_j$ with respect to $\mc'$.
We thus have $t_{ij} = \delta_{ij}$
for all $j\neq k$.
There is a triangle
unique up to isomorphism
$$
M_k' \gfl B_{M_k} \gfl M_k \gfl \susp M_k'
$$
where $B_{M_k} \gfl M_k$ is a right
$\tc\cap\tc'$-approximation.
Since the Auslander--Reiten quiver of
$\tc$ has no loops and no 2-cycles,
$B_{M_k}$ is isomorphic to the direct sum:
$\bigoplus_{i\in I}(M_j')^{a_{i k}}$.
We thus have $t_{ik} = -\delta_{i k} + a_{ik}$,
which equals $\frac{|b_{i k}| + b_{i k}}{2}$.
Remark that, by lemma 7.1 of~\cite{GLSrigid}, as stated 
in section~\ref{subsection: FZ}, 
we have $T^2 = Id$, so that $S = T^{\text{t}}$ and
$$
s_{ij} = \left\{ \begin{array}{ll}
 -\delta_{i j} + \frac{|b_{ij}| - b_{ij}}{2} & \text{if } i=k \\
 \;\;\delta_{ij} & \text{else.}
\end{array} \right.
$$
\end{proof}

\begin{theo}
 The matrix $B_{\mc'}$ is obtained from the matrix 
$B_{\mc}$ by the Fomin--Zelevinski mutation rule
in the direction $M$.
\end{theo}

\begin{proof}
 By~\cite{BFZ} (see section~\ref{subsection: FZ}),
and by lemma~\ref{lem: S}, 
we know that the mutation of the matrix $B_{\mc}$ in 
direction $M$ is given by $T B_{\mc'}T^\text{t}$, 
which is $B_{\mc}$, by the generalized
mutation rule (theorem~\ref{theo: gmr}).
\end{proof}

\subsection{Cluster categories}

In~\cite{BKL}, the authors study the Grothendieck 
group of the cluster category $\cat_A$ associated 
to an algebra $A$ which is either hereditary or 
canonical, endowed with any admissible triangulated 
structure. A triangulated structure on the category 
$\cat_A$ is called admissible in~\cite{BKL} if the 
projection functor from the bounded derived category 
$\dba$ to $\cat_A$ is exact (triangulated). They define 
a Grothendieck group $\kred$ 
with respect to the triangles induced by those of $\dba$, and 
show that it coincides with the usual Grothendieck group of 
the cluster category in many cases: 

\begin{theo}\label{theo: BKL}\emph{[Barot--Kussin--Lenzing]} 
 We have $\kzero(\cat_A) = \kred$ in each of the
following three cases:
\begin{itemize}
\item[(i)] $A$ is canonical with weight sequence $(p_1,\ldots,p_t)$
  having at least one even weight.
\item[(ii)] $A$ is tubular,
\item[(iii)] $A$ is hereditary of finite representation type.
\end{itemize}
\end{theo}

Under some restriction on the triangulated structure 
of $\cat_A$, we have the following generalization of 
case (iii) of theorem~\ref{theo: BKL}:

\begin{theo}
 Let $A$ be a finite-dimensional hereditary algebra, 
and let $\cat_A$ be the associated cluster category 
with its triangulated structure defined in~\cite{Ktri}. 
Then we have $\kzero(\cat_A) = \kred$.
\end{theo}

\begin{proof}
By lemma 3.2 in~\cite{BKL}, this
theorem is a corollary of
the following lemma.
\end{proof}

\begin{lem}
 Under the assumptions of section~\ref{subsection: gmr}, 
and if moreover $\ms$ has a finite number 
$n$ of non-isomorphic indecomposable 
objects, then we have an isomorphism $\kzero(\cat) \simeq 
\zb^n/\im B_{\mc}$.
\end{lem}

\begin{proof}
 This is a restatement of
theorem~\ref{theo: kzero}.
\end{proof}

\nocite{*}
\bibliographystyle{plain}
\bibliography{GeneralizedMutationRule}

\end{document}